\documentclass{amsart}
\usepackage{amssymb,amsmath}
\newcommand{\Real}{\mathbb{R}}

\newcommand{\essinf}{\mathop{\mathrm{ess\;\!inf}}}
\newcommand{\ds}{\displaystyle}
\newcommand{\ie}{\emph{i.e.}}

\newcommand{\cf}{\emph{cf}}
\newcommand{\dist}{\mathrm{dist}}
\newcommand{\vol}{\mathrm{vol}}
\newtheorem{Theorem}{Theorem}

\newtheorem{Proposition}{Proposition}
\theoremstyle{definition}
\newtheorem{Remark}{Remark}
\newtheorem{examp}{Example}
\begin{document}
%
\title[The first Dirichlet eigenvalue and the isoperimetric constant]{\textbf{
A sharp upper bound for the first Dirichlet eigenvalue
and the growth of the isoperimetric constant of convex domains
}}
\author{Pedro Freitas \ and \ David Krej\v{c}i\v{r}{\'\i}k}

\thanks{
Partially supported by FCT, Portugal, through programs
POCTI/MAT/60863/2004, POCTI/POCI2010 and SFRH/BPD/11457/2002.
The second author was also supported by
the Czech Academy of Sciences and its Grant Agency
within the projects IRP AV0Z10480505 and A100480501,
and by the project LC06002 of the Ministry of Education,
Youth and Sports of the Czech Republic.
}

\address{
Department of Mathematics,
Faculdade de Motricidade Humana (TU Lisbon)
{\rm and}
Group of Mathematical Physics of the University of Lisbon,
Complexo Interdisciplinar,
Av.~Prof.~Gama Pinto~2, P-1649-003 Lisboa,
Portugal
}
\email{freitas@cii.fc.ul.pt}
\address{
Department of Theoretical Physics,
Nuclear Physics Institute,
Academy of Sciences,
250\,68 \v{R}e\v{z}, Czech Republic
}
\email{krejcirik@ujf.cas.cz}
\date{October 29, 2007}
%
\begin{abstract}
We show that as the ratio between the first Dirichlet eigenvalues
of a convex domain and of the ball with the same volume becomes
large, the same must happen to the corresponding
ratio of isoperimetric constants.
The proof is based on the generalization to arbitrary dimensions of
P\'{o}lya and Szeg\"{o}'s $1951$ upper bound for the
first eigenvalue of the Dirichlet Laplacian on planar
star-shaped domains which depends on the support function of the domain.

As a by-product, we also obtain a sharp upper bound for the
spectral gap of convex domains.
\end{abstract}
\maketitle
%
%
%
\section{Introduction}
%
Let~$\Omega$ be a bounded domain (\ie~open connected set)
in the $d$-dimensional Euclidean space~$\Real^d$, with $d \geq 1$,
and denote by $\lambda_{1}(\Omega)$ the first eigenvalue of the
Dirichlet Laplacian in $L^2(\Omega)$. The Faber-Krahn inequality
provides a lower bound for $\lambda_{1}(\Omega)$, namely,
\begin{equation}\label{FK}
  \lambda_{1}(\Omega)
  \ \geq \
  \lambda_{1}(B_1) \left(\frac{|B_1|}{|\Omega|}\right)^{2/d},
\end{equation}
where~$B_1$ is the $d$-dimensional ball of unit radius
and the absolute-value
signs denote $d$-dimensional Lebesgue measure
(later on we use the same notation for $(d-1)$-dimensional
Hausdorff measure of the boundary of~$\Omega$ as well).
Modulus a set of zero capacity, equality in~(\ref{FK}) is attained if
and only if $\Omega$~is a ball and it is thus of interest to understand
how strong this connection is. More precisely, and assuming for the time
being that we have a way of measuring how far a set is from the ball of the
same volume in terms of elementary geometric quantities, we would like
to answer questions such as the following:
if a set is {\it far away} from
the ball, must its first Dirichlet eigenvalue be much larger than that
of the ball of the same volume?

This particular question was given a
positive answer in~\cite{gugg}, where the measure of deviation of a
convex domain from the ball which was used was based on the support
function of the domain.
More recently, in~\cite{mela,beco,FMP} the authors addressed
the question of whether a domain for which $\lambda_{1}(\Omega)|\Omega|^{2/d}$
is close to the corresponding quantity for the ball must be close to
the ball in the sense of Fraenkel asymmetry
(\ie\ Hausdorff distance if $\Omega$ is convex),
again providing a positive answer.

In this paper, we consider the issue
of whether having a large first Dirichlet eigenvalue
implies being away from the corresponding ball. Using the trivial upper bound
$\lambda_{1}(\Omega)\leq\lambda_{1}(B_{\rho_{\Omega}})$,
where $\rho_\Omega$ is the inradius of~$\Omega$, it is possible to give
an immediate answer to this question in terms of $\rho_{\Omega}$. However, this
bound is not very good in general, say for long parallelepipeds,
and so our purpose was to obtain a different characterization which would behave
better precisely when away from the ball.

The main result of this paper in this direction
is the following estimate using the isoperimetric constant
as a measure of deviation of~$\Omega$ from~$B$:
\begin{Theorem}\label{thm.isoperimetric}
Let $\Omega$ be a bounded convex domain of $\Real^{d}$.
Then
$$
  \frac{|\partial\Omega|}{|\Omega|^{1-1/d}}
  \ \geq \
  \frac{|\partial B|}{|B|^{1-1/d}} \,
  \sqrt{\frac{\lambda_1(\Omega)}{\lambda_1(B)}} \,
  \frac{\pi}{2\sqrt{\lambda_1(B_1)}}
  \,,
$$
where~$B$ is the ball of volume $|\Omega|$.
\end{Theorem}

The proof is based on the following result:
\begin{Theorem}\label{Thm.Corol}
Let~$\Omega$ be a bounded convex domain of~$\Real^d$.
Then
\[
  \lambda_{1}(\Omega)
  \ \leq \
  \lambda_{1}(B_1) \, \frac{|\partial\Omega|}{d\,\rho_\Omega\,|\Omega|}
  \,.
\]
\end{Theorem}

This upper bound is a consequence
of a stronger upper bound for $\lambda_1(\Omega)$
holding in the more general case of star-shaped domains and which we
believe to be of interest in its own right
(\cf~Theorem~\ref{Thm.main} below). This is an extension
to arbitrary dimensions of an upper bound for $\lambda_1(\Omega)$
appearing in P\'{o}lya and Szeg\"{o}'s $1951$ book~\cite{posz} in the planar case.
As in the case of~\cite{gugg}, this bound also depends on the support function
of the domain in a non-elementary way. Due to this, we postpone the statement of
this result to the next section where we provide the necessary background,
its proof being then given in Section~\ref{proofs}.

The proof of Theorem~\ref{thm.isoperimetric} together with a brief discussion
of optimality, and other applications of our bounds
are given in Section~\ref{Sec.appl} where, in particular, we obtain a sharp upper
bound for the spectral gap. Finally, in the last section we recall some two-dimensional
upper bounds and conjectures which will be used for comparison and consider some examples.

\section{An upper bound for the first eigenvalue of star-shaped domains}
\label{Sec.main}
%
To state P\'{o}lya and Szeg\"{o}'s result
and its generalization to arbitrary
dimension we need to introduce a geometric quantity
which measures how far away we are from the ball,
and which may be expressed in terms of the support function
of the given domain.
To be more precise, let $\Omega$~be a \emph{star-shaped} domain
with respect to a point $\xi\in\Omega$,
\ie, for each point $x\in\partial\Omega$ the segment joining~$\xi$
with~$x$ lies in $\Omega\cup\{x\}$
and is transversal to~$\partial\Omega$ at the point~$x$.
Assume now that the boundary~$\partial\Omega$ is locally Lipschitz.
Then the outward unit normal vector field
$
  N: \partial\Omega \to \Real^d
$
can be uniquely defined almost everywhere on $\partial\Omega$.
At those points $x\in\partial\Omega$ for which $N(x)$
is uniquely defined, we introduce the \emph{support function}
\begin{equation*}
  h_\xi(x) := (x-\xi) \cdot N(x) \,,
\end{equation*}
where the dot denotes the standard scalar product in~$\Real^d$.
We say that~$\Omega$ is \emph{strictly star-shaped}
with respect to the point $\xi\in\Omega$
if~$\Omega$ is star-shaped with respect to~$\xi$
and the support function is uniformly positive,
\ie,
\begin{equation}\label{Ass.strict}
  \essinf_{x\in\partial\Omega} h_\xi(x) > 0 \,.
\end{equation}
In this case, we shall denote by~$\omega$ the set of points
with respect to which~$\Omega$ is strictly star-shaped,
and define the following intrinsic quantity of the domain
\[
 F(\Omega) := \inf_{\xi\in\omega} \int_{\partial\Omega} h_\xi^{-1}.
\]

Our main result reads as follows:
\begin{Theorem}\label{Thm.main}
Let~$\Omega$ be a bounded strictly star-shaped domain in~$\Real^d$
with locally Lipschitz boundary~$\partial\Omega$.
Then
$$
  \lambda_{1}(\Omega)
  \ \leq \
  \lambda_{1}(B_1) \, \frac{F(\Omega)}{d\,|\Omega|}
  \,.
$$
\end{Theorem}
\begin{Remark}
If $d=2$ Theorem~\ref{Thm.main} coincides
with the P\'{o}lya-Szeg\"{o} bound~\cite[Sec.~5.6]{posz}.
\end{Remark}
\begin{Remark}
Combining the upper bound of Theorem~\ref{Thm.main}
with the Faber-Krahn inequality given in~(\ref{FK}),
we see that $F(\Omega)$ is bounded from below by
\[
  F(\Omega)
  \geq
  \frac{\ds d \, |B_{1}|^{2/d}}{\ds |\Omega|^{2/d-1}}
  \,,
\]
with equality only when~$\Omega$ is a ball.
In the two-dimensional case
this was shown in~\cite{aiss} by a different method.
\end{Remark}

If we now restrict ourselves to convex domains
it is possible to simplify the discussion somewhat.
To begin with, we have that the boundary of a convex open
subset of~$\Real^d$ is locally Lipschitz
(\cf~\cite[Sec.~V.4.1]{Edmunds-Evans}).
Furthermore, $\omega=\Omega$. Indeed, for any $\xi\in\Omega$ one has
\begin{equation}\label{distance}
  \essinf_{x\in\partial\Omega} h_\xi(x)
  \geq \dist(\xi,\partial\Omega)
  \qquad
  \mbox{if $\Omega$ is convex},
\end{equation}
which follows from the geometrical meaning of~$h_\xi(x)$
being the distance from~$\xi$
to the tangent space $T_x(\partial\Omega)$.
Finally, (\ref{distance})~can be used to obtain
a simple upper bound to $F(\Omega)$
\begin{equation}\label{inradius}
  F(\Omega) \leq \frac{|\partial\Omega|}{\rho_\Omega}
  \qquad
  \mbox{if $\Omega$ is convex},
\end{equation}
where~$\rho_\Omega$ is the inradius of~$\Omega$.
This approximation readily establishes Theorem~\ref{Thm.Corol},
a weaker version of Theorem~\ref{Thm.main}
but, on the other hand, it allows us to write
the upper bound explicitly in
terms of more elementary geometric quantities.

\begin{Remark}
Since~\cite[Thm.~35.1.2]{buza}
\begin{equation}\label{rhoineq}
  \rho_{\Omega} \, |\partial\Omega|\leq d \, |\Omega|
  \,,
\end{equation}
we see that Theorem~\ref{Thm.Corol} is, in general,
an improvement to the trivial upper bound
given by $\lambda_{1}(\Omega)\leq\lambda_{1}(B_{\rho_{\Omega}})$.
We note, however, that there are certain classes of domains
such as triangles and regular polygons
for which equality holds in~(\ref{rhoineq})
(for a recent discussion of these classes
in the two- and three- dimensional cases
we refer to~\cite{apmn2,apmn1}, respectively).
Moreover, for triangles it can be checked directly
that we also have equality in~(\ref{inradius}),
so that Theorem~\ref{Thm.main} is equivalent
to Theorem~\ref{Thm.Corol}
and our method gives just the trivial result in this case
(see~\cite{aiss} for explicit formulae).
\end{Remark}
\begin{Remark}
When the domain is star-shaped but not necessarily convex,
the bound given by Theorem~\ref{Thm.main} may be worse
than the trivial bound,
as exemplified by the Swiss Cross example considered in~\cite{hesa}
-- see Example~\ref{swcr} below.
\end{Remark}

Although our proof of Theorem~\ref{Thm.Corol}
only holds for the convex case,
we conjecture that this bound holds
for general bounded domains in $\Real^{d}$.

\section{Proof of Theorem~\ref{Thm.main}\label{proofs}}
%
Henceforth let us assume that $d \geq 2$
(for one-dimensional intervals
we clearly have equalities in the upper bounds
of Theorems~\ref{Thm.Corol} and~\ref{Thm.main}).

We need to introduce some notation concerning
the geometry of the hypersurface~$\partial\Omega$.
By assumption,
each point~$x$ on the boundary $\partial\Omega$
has a neighbourhood (an open subset of~$\Real^d$)
whose intersection with the boundary,
denoted by $V \subset \partial\Omega$,
is $C^{0,1}$-diffeomorphic to an open subset
$U\subset\Real^{d-1}$ by means of a chart $\Gamma:U \to V$.
Having in mind that, by Rademacher's theorem,
$\Gamma$~is differentiable almost everywhere in~$U$,
let $g_{ij}$ denote the coefficients
of the metric tensor of~$\partial\Omega$
induced by these local diffeomorphisms, \ie,
$$
  g_{ij} := (\partial_i\Gamma) \cdot (\partial_j\Gamma)
  \,, \qquad
  i,j\in\{1,\dots,d-1\}
  \,.
$$
Recall that the cross-product
$
  [\partial_1\Gamma,\dots,\partial_{d-1}\Gamma]
$
is perpendicular to~$\partial\Omega$
and that its magnitude is equal to the square root of
$
  |g| := \det(g_{ij})
$.

Let~$\Omega$ be strictly star-shaped with respect
to $\xi\in\Omega$.
We parameterize $\Omega\setminus\{\xi\}$
by means of the mapping
\begin{equation}\label{coordinates}
  \mathcal{L}: \partial\Omega\times(0,1) \to \Real^d:
  \big\{ (x,t) \mapsto \xi + (x-\xi) t \big\}
  \,.
\end{equation}
Notice that the ``shrunk boundary''
$
  \mathcal{L}(\partial\Omega\times\{t\})
$
is indeed contained in~$\Omega$ for any $t\in(0,1)$.
Using the local parametrization of the boundary by~$\Gamma$
and properties of the cross-product,
the Jacobian of the transformation~$\mathcal{L}$
can be locally identified with
$$
  J(\cdot,t) =
  \begin{vmatrix}
  (\partial_1\Gamma^1) \, t & \dots
  & (\partial_{d-1}\Gamma^1) \, t & \Gamma^1
  \\
  \vdots & & \vdots & \vdots
  \\
  (\partial_1\Gamma^d) \, t& \dots
  & (\partial_{d-1}\Gamma^d) \, t & \Gamma^d
  \end{vmatrix}
  = [\partial_1\Gamma,\dots,\partial_{d-1}\Gamma] \cdot \Gamma \ t^{d-1}
  \,.
$$
Hence
\begin{equation}\label{Jacobian}
  |J(u,t)| = |g|^{1/2}(u) \ h_\xi\big(\Gamma(u)\big) \ t^{d-1}
\end{equation}
for every $t\in(0,1)$ and almost every $u \in U$.
By virtue of the inverse function theorem
and assumption~(\ref{Ass.strict}),
we therefore conclude that
$
  \mathcal{L} : \partial\Omega \times (0,1)
  \to \Omega\setminus\{\xi\}
$
is indeed a diffeomorphism.

In other words, the Euclidean domain~$\Omega$
with the point~$\xi$ removed
can be identified with the Riemannian manifold
$$
  M:=\big(\partial\Omega \times (0,1),G\big) \,,
$$
where the metric~$G$ is induced by~(\ref{coordinates}).
Since the coefficients of~$G$ are locally given by
$$
  G_{ij} := (\partial_i\mathfrak{L}) \cdot (\partial_j\mathfrak{L})
  \,, \qquad
  i,j\in\{1,\dots,d\}
  \,,
$$
where
$
  \mathfrak{L}:=\mathcal{L}\circ(\Gamma \otimes 1)
$
with~$1$ being the identity function on the interval $(0,1)$,
we plainly have
\begin{equation*}
  \big(G_{ij}(\cdot,t)\big) =
  \begin{pmatrix}
  g_{11} \, t^2 & \dots & g_{1 d-1} \, t^2
  & \Gamma\cdot(\partial_1\Gamma) \, t
  \\
  \vdots & & \vdots & \vdots
  \\
  g_{d-1 1} \, t^2 & \dots & g_{d-1 d-1} \, t^2
  & \Gamma\cdot(\partial_{d-1}\Gamma) \, t
  \\
  \Gamma\cdot(\partial_1\Gamma) \, t
  & \dots & \Gamma\cdot(\partial_{d-1}\Gamma) \, t
  & |\Gamma|^2
  \end{pmatrix}
  \,.
\end{equation*}
Using the relation
$
  |G| := \det(G_{ij}) = J^2
$
and formula~(\ref{Jacobian}),
we see that the volume element $d\vol$ of the manifold~$M$
is decoupled as follows:
\begin{equation}\label{volume}
  d\vol = h_\xi(x) \, d\sigma(x) \ t^{d-1} \, dt
  \,,
\end{equation}
where $d\sigma$ and $dt$ denote the respective
measures on $\partial\Omega$ and $(0,1)$.
At the same time,
denoting by $G^{ij}$ the coefficients
of the matrix inverse to~$(G_{ij})$,
we locally find that
$$
  G^{dd}(u,t) = h_\xi^{-2}\big(\Gamma(u)\big)
  \,.
$$
If~$\psi$ is any differentiable function on $(0,1)$
and~$1$ denotes the identity function on~$\partial\Omega$,
then the last result implies
\begin{equation}\label{gradient}
  \big|\nabla_{\!G}(1\otimes\psi)\big|_G = h_\xi^{-1} \, |\psi'|
  \,,
\end{equation}
where~$\nabla_{\!G}$ and~$|\cdot|_G$ stand for
the gradient and norm on~$M$, respectively.

Using the above geometric preliminaries,
the Hilbert space~$L^2(\Omega)$
can be identified with
$
  L^2(M) \equiv
  L^2\big(\partial\Omega\times(0,1),d\vol\big)
$
and the Dirichlet Laplacian in the former
is unitarily equivalent to the self-adjoint operator
associated in the latter with the quadratic form
$$
  Q[\Psi] :=
  \big\| \,|\nabla_{\!G}\Psi|_G \big\|_{L^2(M)}^2
  \,, \qquad
  \Psi \in D(Q) := H_0^1(M)
  \,.
$$
Recall that the Sobolev space $H_0^1(M)$ is the completion of
$C_0^\infty\big(\partial\Omega\times(0,1)\big)$
with respect to the norm
$
  ( Q[\cdot]+\|\cdot\|_{L^2(M)}^2 )^{1/2}
$.
Let~$\psi$ be a non-zero function from
$H_0^1\big((0,1),t^{d-1}\,dt\big)$
and let~$1$ denote the identity function on~$\partial\Omega$.
Then the function $1\otimes\psi$
belongs to the form domain $D(Q)$
and we can use it as a test function
in the variational formulation
for the first eigenvalue
of the operator associated to~$Q$.
Employing~(\ref{volume}) and~(\ref{gradient}),
we get
\begin{equation*}
  \lambda_{1}(\Omega)
  \ \leq \
  \frac{Q[1\otimes\psi]}{\ \|1\otimes\psi\|_{L^2(\Omega)}^2}
  \ = \
  \frac{\int_{\partial\Omega} h_\xi(x)^{-1} d\sigma(x)}
  {\int_{\partial\Omega} h_\xi(x) \, d\sigma(x)}
  \
  \frac{\int_0^1 |\psi'(t)| \, t^{d-1} dt}
  {\int_0^1 |\psi(t)|^2 \, t^{d-1} dt}
  =: \lambda_{1}(\Omega;\psi,\xi)
  \,.
\end{equation*}

Now, if~$\Omega$ is the ball of radius~$1$ centered at~$\xi$,
the above result reduces to
$$
  \lambda_{1}(B_1) \leq \frac{\int_0^1 |\psi'(t)| \, t^{d-1} dt}
  {\int_0^1 |\psi(t)|^2 \, t^{d-1} dt}
  \,.
$$
But we know that the equality sign holds
if, and only if, $\psi$~is chosen as the radial
component of the first eigenfunction
of the Dirichlet Laplacian in the ball.
Indeed, this eigenfunction is radially symmetric
and can be written as $1\otimes\psi$ in our coordinates.
Consequently,
$$
  \min_\psi \lambda_{1}(\Omega;\psi,\xi)
  = \lambda_{1}(B_1) \
  \frac{\int_{\partial\Omega} h_\xi(x)^{-1} d\sigma(x)}
  {\int_{\partial\Omega} h_\xi(x) \, d\sigma(x)}
  \,.
$$
Minimizing the integral in the numerator with respect to~$\xi$,
we arrive at the quantity~$F(\Omega)$ of Theorem~\ref{Thm.main}.
It remains to realize that the integral of support function
is actually independent of~$\xi$
because, by~(\ref{volume}) and Fubini's theorem,
$$
  |\Omega| =
  \int_{\partial\Omega\times(0,1)} d\vol
  = \frac{1}{d} \, \int_{\partial\Omega} h_\xi(x) \, d\sigma(x)
  \,.
$$
This concludes the proof of Theorem~\ref{Thm.main}.

\section{Applications}\label{Sec.appl}
%
\subsection{Torsional rigidity}
The method of the present paper applies
to other Sobolev-inequality-type problems, too.
For instance,
let $P(\Omega)$ be the \emph{torsional rigidity} of~$\Omega$
defined by
$$
  \frac{1}{P(\Omega)}
  := \inf_{\psi \in H_0^1(\Omega)\setminus\{0\}}
  \frac{\|\nabla\psi\|_{L^2(\Omega)}^2}
  {\|\psi\|_{L^1(\Omega)}^2}
  \,.
$$
Then, following the lines of Section~\ref{proofs}, we have
\begin{Proposition}
Let~$\Omega$ be a bounded strictly star-shaped domain in~$\Real^d$
with locally Lipschitz boundary~$\partial\Omega$.
Then
$$
  \frac{1}{P(\Omega)}
  \ \leq \
  \frac{|B_1|}{P(B_1)} \, \frac{F(\Omega)}{d\,|\Omega|^2}
  \,.
$$
\end{Proposition}
\noindent
Again, the equality is attained for~$\Omega$ being a ball.
For $d=2$ this result coincides
with the P\'olya-Szeg\"o bound \cite[Sec.~5.5]{posz}.

\subsection{The growth of the isoperimetric constant}
We shall now prove Theorem~\ref{thm.isoperimetric}.
The idea is to estimate the inradius appearing
in the bound of Theorem~\ref{Thm.Corol} by means of
the following lower bound due to Protter~\cite{prot}
\[
\lambda_{1}(\Omega)\geq \frac{\ds \pi^{2}}{\ds 4\rho^{2}_{\Omega}}
\,.
\]
(Protter's bound actually includes an extra term
depending on the diameter, but for our purposes
it is sufficient to consider the expression above.)
This leads to
$$
  \frac{|\partial\Omega|}{d\,|\Omega|}
  \geq
  \sqrt{\frac{\lambda_1(\Omega)}{\lambda_1(B_1)}} \,
  \frac{\pi}{2\sqrt{\lambda_1(B_1)}}
  \,,
$$
which is equivalent to the inequality of Theorem~\ref{thm.isoperimetric}
due to the scaling properties
$\lambda_1(B)=\lambda_1(B_1) r^{-2}$ and $|B|=|B_1|r^d$
where~$r$ is the radius of~$B$,
and $|\partial B_1|=d |B_1|$.

Since Protter's bound used in the proof is not sharp for the ball,
the inequality in the above theorem is not an improvement
upon the classical isoperimetric inequality for $\lambda_{1}(\Omega)/\lambda_1(B)$
close to one.
However, it is clear that this will be the case when this ratio becomes large.
At the same time, our bound is optimal in the sense that
there exist $d$-dimensional domains for which the growth
of $\sqrt{\lambda_1(\Omega)/\lambda_1(B)}$ cannot be improved.
More precisely, it is not difficult to obtain
that if~$\mathcal{R}$ is a $d$-dimensional parallelepiped
we have
$$
  \frac{|\partial\mathcal{R}|}{|\mathcal{R}|^{1-1/d}}
  \leq
  \frac{|\partial B|}{|B|^{1-1/d}} \,
  \sqrt{\frac{\lambda_1(\mathcal{R})}{\lambda_1(B)}} \,
  \frac{2\sqrt{\lambda_1(B_1)}}{\pi\sqrt{d}}
  \,.
$$
The question of the optimal constant multiplying $\sqrt{\lambda_1(\Omega)/\lambda_1(B)}$
in Theorem~\ref{thm.isoperimetric}
remains open.

\subsection{The second eigenvalue and the gap\label{sec:gap}}
Combining the bound in Theorem~\ref{Thm.Corol} with the Ashbaugh-Benguria
bound~\cite{asbe} for the spectral
quotient $\lambda_{2}(\Omega)/\lambda_{1}(\Omega)$ gives a similar upper
bound for the second eigenvalue of a convex domain.
This, together with the Faber-Krahn inequality
yields, in turn, an upper bound for the spectral gap.
\begin{Proposition}
Let $\Omega$ be a bounded convex domain of $\Real^{d}$.
Then the second Dirichlet eigenvalue $\lambda_{2}(\Omega)$ satisfies
\[
\lambda_{2}(\Omega)
  \ \leq \
  \lambda_{2}(B_1) \, \frac{|\partial\Omega|}{d\,\rho_\Omega\,|\Omega|}
  \,.
\]
As a consequence, the spectral gap satisfies
\[
\lambda_{2}(\Omega)-\lambda_{1}(\Omega)
\ \leq \
\lambda_{2}(B_1) \, \frac{|\partial\Omega|}{d\,\rho_\Omega\,|\Omega|} -
\lambda_{1}(B_1) \left(\frac{|B_1|}{|\Omega|}\right)^{2/d}.
\]
\end{Proposition}
For simplicity we used the bound for convex sets, but of
course that the bound in Theorem~\ref{Thm.main} provides a better result
valid for star-shaped domains.
Note also that both upper bounds give equality for the ball.
For a numerical study regarding upper and lower bounds of the gap
see~\cite{anfr2}, which also contains some new conjectures for this
problem.

\section{Discussion of the upper bounds and examples}\label{examples}
%
In spite of the fact that due to the existence
of Rayleigh's variational formulation
upper bounds for the first eigenvalue are, in principle,
easier to obtain than lower bounds,
it is slightly more delicate to obtain sharp upper bounds
depending on the geometric quantities used here
and which are valid in arbitrary dimensions.
We shall thus now compare Theorems~\ref{Thm.Corol} and~\ref{Thm.main} to
other existing bounds and conjectures for some particular examples.

With this in mind we begin by recalling some two-dimensional bounds.
Among these there is the family of upper bounds based on the
method of parallel coordinates,
which includes P\'{o}lya's $1960$ bound
for simply-connected domains~\cite{poly}
(sharp asymptotically on infinite rectangular strips)
\[
  \lambda_{1}(\Omega)
  \ \leq \
  \frac{\ds \pi^{2}}{\ds 4}
  \frac{\ds |\partial \Omega|^{2}}{\ds |\Omega|^{2}}
  \,,
\]
and its improvement by Payne and Weinberger~\cite{pawe}, namely,
\begin{equation}\label{pwbound}
  \lambda_{1}(\Omega)
  \ \leq \
  \frac{\ds \pi j_{0,1}^{2}}{\ds |\Omega|}\left[
  1+\left(\frac{\ds 1}{\ds
  J_{1}^{2}(j_{0,1})}-1\right)\left(\frac{\ds |\partial\Omega|^2}{\ds 4\pi
  |\Omega|}-1\right)\right].
\end{equation}
Here $J_{1}$ and $j_{0,1}$ denote, respectively, the Bessel function of
the first kind of order one, and the first positive zero of $J_{0}$,
the Bessel function of the first kind of order zero. The Payne-Weinberger
bound~(\ref{pwbound}) is an explicit expression in terms of the area and
perimeter which is obtained from their stronger bound
\begin{equation}\label{PW-bound}
  \lambda_{1}(\Omega)
  \ \leq \
  \frac{4\pi^2}{|\partial\Omega|^2} \ k(p)^2
  \,, \qquad
  p := 1-\frac{4\pi|\Omega|}{|\partial\Omega|^2}
  \,,
\end{equation}
where~$k=k(p)$ is the first zero of the transcendental equation
$$
  J_0(k) Y_1(pk) = Y_0(k) J_1(pk)
  \,.
$$
Here $Y_{0}$ and $Y_1$ denote the Bessel functions
of the second kind of order zero and one, respectively.
It is this stronger bound that we will consider throughout this section
for comparison,
and we shall refer to it as the PW-bound.
Note that~(\ref{PW-bound}) is sharp for the disc
and asymptotically for infinite rectangular strips.
While a generalization of P\'olya's bound to arbitrary dimensions
can be found in~\cite{Savo},
the proof of the stronger result~(\ref{PW-bound}) does not seem
to have a straightforward extension to higher dimensions.

We also remark that although~(\ref{pwbound})
does give equality on the disc,
the numerical study carried out in~\cite{anfr}
suggests that this bound might
still be improved and it is conjectured there
that the optimal bound depending
explicitly on the area and the perimeter and valid
for simply-connected two-dimensional domains should be
\begin{equation}\label{afconj}
  \lambda_{1}(\Omega)
  \ \leq \
  \frac{\ds \pi j_{01}^2}{\ds |\Omega|}
  +\frac{\ds \pi^2}{\ds 4}
  \frac{\ds  |\partial\Omega|^2 - 4\pi |\Omega|}{\ds |\Omega|^2}
  \,,
\end{equation}
providing now equality not only for the disc
but also asymptotically on infinite rectangular strips.

Along different lines, Maz'ya and Shubin have recently proved
upper and lower bounds depending on the
interior capacity radius~\cite{mash}.

We shall now consider some examples for which we compare
the upper bounds given by Theorems~\ref{Thm.Corol} and~\ref{Thm.main}
with the PW-bound~(\ref{PW-bound}) and conjecture~(\ref{afconj}).

\begin{examp}[Rectangular parallelepipeds]
%
Given positive numbers $a_1, \dots, a_d$,
let
$
  \mathcal{R} := (-a_1,a_1)\times\dots\times(-a_d,a_d)
$.
Elementary calculations show that the infimum
in the definition of~$F(\mathcal{R})$ is attained
for the intuitive choice $\xi=0$,
with the result
$$
  F(\mathcal{R})
  = |\mathcal{R}| \left(a_1^{-2}+\dots+a_d^{-2}\right)
  \,.
$$

For rectangles, conjecture~(\ref{afconj})
is better than the PW-bound~(\ref{PW-bound})
for all the values of the parameter $c:=a_1/a_2\in(0,1]$.
Theorem~\ref{Thm.main} (respectively Theorem~\ref{Thm.Corol})
provides a better upper bound than conjecture~(\ref{afconj})
in the range of $c\in(0.3,1]$ (respectively $c\in(0.7,1]$).
The largest discrepancy between the upper bound of Theorem~\ref{Thm.main}
(respectively conjecture~(\ref{afconj}))
and the actual value of $\lambda_{1}(\mathcal{R})$
is in the limit $c \to 0$ (respectively for $c=1$)
when it is about $15\%$ (respectively $26\%$).
\end{examp}

\begin{examp}[Ellipsoids]
%
Given positive numbers $a_1, \dots, a_d$,
let~$\mathcal{E}$ be the domain enclosed by an ellipsoid,
\ie~the surface determined by the implicit equation
$
  f(x):=
  (x_1/a_1)^2+\dots+(x_d/a_d)^2-1=0
$.
First of all, by symmetry,
it is possible to conclude that
the infimum in the definition of~$F(\mathcal{E})$
is attained for~$\xi=0$.
Since $\nabla f/|\nabla f|$ is either $+N$ or $-N$
uniformly on~$\partial\Omega$, we have
$$
  h_0^{-1}(x)
  = N(x) \cdot \frac{\nabla f(x)}{x\cdot\nabla f(x)}
  = N(x) \cdot \left(\frac{x_1}{a_1^2},\dots,\frac{x_d}{a_d^2}\right)
  .
$$
Now using the divergence theorem,
we arrive at
$$
  F(\mathcal{E})
  = |\mathcal{E}| \left(a_1^{-2}+\dots+a_d^{-2}\right)
  \,.
$$
That is, we formally obtain the same upper bound
as in the case of parallelepipeds
(notice that the volume terms in the bound of
Theorem~\ref{Thm.main} cancel).
However, the present bound is better because
$\mathcal{E}\subset\mathcal{R}$,
so that
$
  \lambda_{1}(\mathcal{E}) \geq \lambda_{1}(\mathcal{R})
$
by monotonicity of Dirichlet eigenvalues.

For ellipses,
Theorem~\ref{Thm.main} provides a better upper bound
than conjecture~(\ref{afconj})
(which is again better than the PW-bound~(\ref{PW-bound}))
for all the values of the parameter $c:=a_1/a_2\in(0,1]$.
Theorem~\ref{Thm.Corol} provides a better upper bound
than the conjecture~(\ref{afconj}) in the regime of $c\in(0,0.1]$.
In fact, the upper bound of Theorem~\ref{Thm.main}
for two-dimensional ellipses seems to be familiar
in the applied sciences
and it is known that for $c \geq 0.5$
the discrepancy between this
and $\lambda_{1}(\mathcal{E})$ does not exceed~$1\%$
(\cf~\cite[Sec.~7.3.4--3]{encyclopedia}).

Finally, let us mention that we obtain the same
formula for $F(\mathcal{E})$
(and therefore for the upper bound of Theorem~\ref{Thm.main})
also in the case when~$\mathcal{E}$
is a tube of elliptical cross-section,
\ie~the domain determined by
$f(x_1,\dots,x_{d-1},0)=0$ and $x_d\in(-a_d,a_d)$.
\end{examp}
\begin{examp}[Stadium]
%
Given positive numbers~$a$ and~$b$,
let~$\mathcal{S}\subset\Real^2$
be the union of the rectangle $(-b,b)\times(-a,a)$
and two discs of radius~$a$ centered
at the points $(-b,0)$ and $(b,0)$.
We put $c:=b/a \in [0,+\infty)$.
Again, by symmetry,
it is possible to conclude that
the infimum in the definition of~$F(\mathcal{S})$
is attained for~$\xi=0$.
Since the boundary of~$\mathcal{S}$ is composed
of straight and arc segments,
the integral of the inverse of the support function
can be computed explicitly:
$$
  F(\mathcal{S}) =
  \begin{cases}
    {\displaystyle
    4c +
    \frac{8}{\sqrt{1-c^2}} \, \arctan\sqrt{\frac{1-c}{1+c}}
    }
    & \mbox{if}\quad  c<1 \,, \\
    8
    & \mbox{if}\quad  c=1 \,, \\
    {\displaystyle
    4c +
    \frac{4}{\sqrt{c^2-1}} \, \log\left(c+\sqrt{c^2-1}\right)
    }
    & \mbox{if}\quad  c>1 \,. \\
  \end{cases}
$$
In this example, Theorem~\ref{Thm.main} provides a better upper bound
than conjecture~(\ref{afconj})
(which is again better than the PW-bound~(\ref{PW-bound}))
for $c\in[0,2.6]$, while Theorem~\ref{Thm.Corol} is worse than
both the conjecture and PW-bound for all the values of the parameter.

It is also possible to consider the asymmetric domain
$\{(x_1,x_2)\in\mathcal{S}\,|\,x_1>0\}$.
Then the position of~$\xi$ minimizing
the infimum in the definition of~$F(\mathcal{S})$
significantly depends on the value of~$c$.
\end{examp}

\begin{examp}[Swiss cross]\label{swcr}
%
As an example of a non-convex domain
(but strictly star-shaped with respect to the origin),
let~$\mathcal{C}\subset\Real^2$
be the union of the two rectangles
$(-b-a,b+a)\times(-a,a)$ and $(-a,a)\times(-b-a,b+a)$.
We put $c:=b/a \in [0,+\infty)$.
An explicit calculation yields
$$
  F(\mathcal{C})= 8 \, \frac{1+c+c^2}{1+c}
  \,.
$$
In this example, Theorem~\ref{Thm.main} provides a better upper bound
than conjecture~(\ref{afconj})
(which is again better than the PW-bound~(\ref{PW-bound}))
for $c\in[0,3.8]$.
The case $c=2$ was numerically analysed in~\cite{hesa}
and it was shown that the discrepancy between the bound
and $\lambda_1(\mathcal{C})$ is less than $39\%$.
\end{examp}

%
%
\providecommand{\bysame}{\leavevmode\hbox to3em{\hrulefill}\thinspace}
\providecommand{\MR}{\relax\ifhmode\unskip\space\fi MR }
\providecommand{\MRhref}[2]{%
  \href{http://www.ams.org/mathscinet-getitem?mr=#1}{#2}
}
\providecommand{\href}[2]{#2}

\end{document}